\documentstyle[12pt]{article}
\newcommand{\sect}[1]{\section{#1}\setcounter{equation}{0}}

\font\mbn=msbm10 scaled \magstep1
\font\mbs=msbm7 scaled \magstep1
\font\mbss=msbm5 scaled \magstep1
\newfam\mbff
\textfont\mbff=\mbn
\scriptfont\mbff=\mbs
\scriptscriptfont\mbff=\mbss\def\mbf{\fam\mbff}
\def\Re{{\mbf R}}

\def\Co{{\mbf C}}

\newtheorem{Th}{Theorem}[section]

\newtheorem{D}[Th]{Definition}
\newtheorem{Proposition}[Th]{Proposition}
\newtheorem{R}[Th]{Remark}
\newtheorem{E}[Th]{Example}
\author{Alexander Brudnyi\thanks{Research supported in part by NSERC.
}\\
Department of Mathematics and Statistics\\
University of Calgary, Calgary\\
Canada}
\title{Holomorphic $L^{p}$-functions on Coverings of Strongly Pseudoconvex Manifolds}
\date{} 
\begin{document} 
\maketitle
\begin{abstract}
{In this paper we will show how to 
construct holomorphic $L^{p}$-functions on unbranched coverings of strongly pseudoconvex manifolds. Also, we prove some extension and approximation theorems for such functions.}
\end{abstract}
\sect{\hspace*{-1em}. Introduction.}
{\bf 1.1.} The present paper continues the study of holomorphic functions of slow growth on unbranched coverings of strongly pseudoconvex manifolds started in [Br1]-[Br3]. Our work was inspired by the seminal paper [GHS] of Gromov, Henkin and Shubin on holomorphic $L^{2}$-functions on coverings of pseudoconvex manifolds. A particular interest in this subject is because of its possible applications to the Shafarevich conjecture on holomorphic convexity of universal coverings of complex projective manifolds. The results of this paper don't imply directly any new results in the area of the Shafarevich conjecture. However, one obtains a rich complex function
theory on coverings of strongly pseudoconvex manifolds that together with some additional methods and ideas would lead to a progress in this conjecture. 

The main result of [Br3] deals with holomorphic $L^{2}$-functions on unbranched coverings of strongly pseudoconvex manifolds. In the present paper we use this result to construct holomorphic $L^{p}$-functions
($p\neq 2$) on such coverings. Also, we prove some extension and approximation theorems for these functions.
In our proofs we exploit some ideas 
based on infinite-dimensional versions of Cartan's A and B theorems 
originally proved by Bungart [B] (see also [L] and references therein for
some generalizations of results of the complex function theory to the
case of Banach-valued holomorphic functions).\\
{\bf 1.2.} To formulate our results we first recall some basic definitions.

Let $M\subset\subset N$ be a domain with smooth boundary 
$bM$ in an $n$-dimensional complex manifold $N$, specifically,
\begin{equation}\label{m1}
M=\{z\in N\ :\ \rho(z)<0\}
\end{equation}
where $\rho$ is a real-valued function of class $C^{2}(\Omega)$ in a
neighbourhood $\Omega$ of the compact set $\overline{M}:=M\cup bM$
such that
\begin{equation}\label{m2}
d\rho(z)\neq 0\ \ \ {\rm for\ all}\ \ \ z\in bM\ .
\end{equation}
Let $z_{1},\dots, z_{n}$ be complex local coordinates in $N$ near $z\in bM$.
Then the tangent space $T_{z}N$ at $z$ is identified with $\Co^{n}$.
By $T_{z}^{c}(bM)\subset T_{z}N$ we denote the complex tangent space to
$bM$ at $z$, i.e.,
\begin{equation}\label{m3}
T_{z}^{c}(bM)=\{w=(w_{1},\dots,w_{n})\in T_{z}(N)\ :\ \sum_{j=1}^{n}
\frac{\partial\rho}{\partial z_{j}}(z)w_{j}=0\}\ .
\end{equation}
The {\em Levi form} of $\rho$ at $z\in bM$ is a hermitian form on
$T_{z}^{c}(bM)$ defined in the local coordinates by the formula
\begin{equation}\label{m4}
L_{z}(w,\overline{w})=\sum_{j,k=1}^{n}
\frac{\partial^{2}\rho}{\partial z_{j}\partial\overline{z}_{k}}(z)w_{j}
\overline{w}_{k}\ .
\end{equation}
The manifold $M$ is called {\em strongly
pseudoconvex} if $L_{z}(w,\overline{w})>0$ for all $z\in bM$ and all
$w\neq 0$, $w\in T_{z}^{c}(bM)$.
Equivalently, strongly pseudoconvex manifolds can be described as the ones 
which locally, in a neighbourhood of any boundary point, can be presented as 
strictly convex domains in $\Co^{n}$. It is also known (see [C], [R]) that 
any strongly pseudoconvex manifold admits a proper holomorphic map with 
connected fibres onto a normal Stein space.\\
{\bf 1.3.}
Without loss of generality we may and will assume that 
$\pi_{1}(M)=\pi_{1}(N)$ for $M$ as above and $N$ is strongly pseudoconvex, as well. (Here $\pi_{1}(X)$ stands for the fundamental group of $X$.)
Let $r: N'\to N$ be an unbranched covering of $N$. For $U\subset N$ we set $U':=r^{-1}(U)$. 

Let $\psi:N'\to\Re_{+}$ be
such that $\log\psi$ is uniformly continuous with respect to the
path metric induced by a Riemannian metric pulled back from $N$.
By ${\cal H}_{p,\psi}(M')$, $1\leq p<\infty$, we denote the
Banach space of holomorphic functions $g$ on $M'$ with norm
\begin{equation}\label{e1.5}
|g|_{p,\psi}:=\sup_{x\in M}\left(\sum_{y\in r^{-1}(x)}
|f(y)|^{p}\psi(y)\right)^{1/p} .
\end{equation}
By ${\cal H}_{\infty,\psi}(M')$ we denote the
Banach space of holomorphic functions $g$ on $M'$ with norm
\begin{equation}\label{e1.6}
|g|_{\infty,\psi}:=\sup_{z\in M'}\{|g(z)|\psi(z)\} .
\end{equation}
\begin{E}\label{ex1}
{\rm Let $d$ be the path metric on $N'$ obtained by
the pullback of a Riemannian metric defined on $N$. Fix
a point $o\in M'$ and set
$$
d_{o}(x):=d(o,x)\ ,\  \ \ x\in N'\ .
$$
It is easy to show by means of the triangle inequality that as the 
function $\psi$ one can take, e.g., $(1+d_{o})^{\alpha}$
or $e^{\alpha d_{o}}$ with $\alpha\in\Re$.}
\end{E}
\begin{R}\label{r1}
{\rm Let $dV_{M'}$ be the Riemannian volume form on the covering $M'$ 
obtained by a Riemannian metric pulled back from $N$. 
Note that every $f\in {\cal H}_{p,\psi}(M')$, $1\leq p<\infty$, also belongs 
to the Banach space $H^{p}_{\psi}(M')$ of holomorphic functions $g$ on
$M'$ with norm
$$
\left(\int_{z\in M'}|g(z)|^{p}\psi(z)dV_{M'}(z)\right)
^{1/p} .
$$
Moreover, one has a continuous embedding 
${\cal H}_{p,\psi}(M')\hookrightarrow H^{p}_{\psi}(M')$.}
\end{R}

Next, we introduce the
Banach space $l_{p,\psi, x}(M')$, $x\in M$, $1\leq p<\infty$,
of functions $g$ on $r^{-1}(x)$ with norm
\begin{equation}\label{e1.7}
|g|_{p,\psi,x}:=\left(\sum_{y\in
r^{-1}(x)}|g(y)|^{p}\psi(y)\right)^{1/p},
\end{equation}
and the Banach space $l_{\infty,\psi, x}(M')$, $x\in M$,
of functions $g$ on $r^{-1}(x)$ with norm
\begin{equation}\label{e1.8}
|g|_{\infty,\psi,x}:=\sup_{y\in
r^{-1}(x)}\{|g(y)|\psi(y)\}.
\end{equation}

Let $C_{M}\subset M$ be the union of
all compact complex subvarieties of $M$ of 
complex dimension $\geq 1$. It is known that if $M$ is strongly
pseudoconvex, then $C_{M}$ is a compact complex subvariety of $M$.

In the sequel for Banach spaces $E$ and $F$, by 
${\cal B}(E,F)$ we denote the space of
all linear bounded operators $E\to F$ with norm $||\cdot||$.

Our main result is the following interpolation theorem.
\begin{Th}\label{te1}
Suppose that $\Omega\subset M\setminus C_{M}$ is an open Stein subset and $K\subset\subset\Omega$. Then for any
$p\in [1,\infty]$ there exists a family 
$\{L_{z}\in {\cal B}(l_{p,\psi,z}(M'),
{\cal H}_{p,\psi}(M'))\}_{z\in\Omega}$ holomorphic 
in $z$ such that 
$$
(L_{z}h)(x)=h(x)\ \ \ {\rm for\ any}\ \ \ h\in l_{p,\psi,z}(M')\ \ \
{\rm and}\ \ \ x\in r^{-1}(z)\ .
$$
Moreover,
$$
\sup_{z\in K}||L_{z}||<\infty\ .
$$
\end{Th}

A similar result for $M$ being a bounded domain in a Stein manifold was proved in [Br4, Theorem 1.3].\\
{\bf 1.4.} To formulate applications of  Theorem \ref{te1} we recall some definitions from [Br4].
\begin{D}\label{holo}
Let $r:N'\to N$ be a covering and $X\subset N$ be a complex submanifold of $N$. 
By ${\cal H}_{p,\psi}(X')$, $X':=r^{-1}(X)$, we denote the Banach space of holomorphic functions $f$ on $X'$ 
such that $f|_{r^{-1}(x)}\in l_{p,\psi,x}(N')$ for any $x\in X$
with norm 
\begin{equation}\label{normx}
\sup_{x\in X}|f|_{r^{-1}(x)}|_{p,\psi,x}\ .
\end{equation}
\end{D}

As an application of Theorem \ref{te1} we prove a result on extension of holomorphic functions from complex submanifolds.

Let $U$ be a relatively compact open subset of a holomorphically convex domain $V\subset\subset N$ containing $C_{M}$ and
$Y\subset V\setminus C_{M}$ be a closed complex 
submanifold of $V$. We set $X:=Y\cap U$. Consider a covering 
$r:N'\to N$. 
\begin{Th}\label{te2}
For every $f\in {\cal H}_{p,\psi}(Y')$, there is
a function $F\in {\cal H}_{p,\psi}(U')$ such that $F=f$ on $X'$.
\end{Th}
\begin{R}\label{re3}
{\rm Let $M\subset\subset N$ be a strongly pseudoconvex manifold. As before, we assume that $\pi_{1}(M)=\pi_{1}(N)$ and
$N$ is strongly pseudoconvex, as well. Then there exist a normal
Stein space $X_{N}$, a proper holomorphic surjective map $p:N\to X_{N}$ with
connected fibres and points $x_{1},\dots, x_{l}\in X_{N}$ such that
$$
p:N\setminus\bigcup_{1\leq i\leq l}p^{-1}(x_{i})\to 
X_{N}\setminus\bigcup_{1\leq i\leq l}\{x_{i}\}
$$ 
is biholomorphic, see [C], [R]. By definition, the domain 
$X_{M}:=p(M)\subset X_{N}$ is strongly pseudoconvex (so it is Stein). Without loss of generality we will assume that  $x_{1},\dots, x_{l}\in X_{M}$ so that
$\cup_{1\leq i\leq l}\ p^{-1}(x_{i})=C_{M}$. Next, $X_{V}:=p(V)$
is a Stein subdomain of $X_{N}$.
Now, as $Y$ we take the preimage under $p$ of a closed complex submanifold of $X_{V}$ that does not contain points $x_{1},\dots, x_{l}$.}
\end{R}

Another application of Theorem \ref{te1} is the following approximation result.
\begin{Th}\label{te3}
Let $K\subset\subset M\setminus C_{M}$ be a compact holomorphically convex subset and $O\subset M\setminus C_{M}$ be a neighbourhood of $K$. Then every function $f\in {\cal H}_{p,\psi}(O')$ can be uniformly  approximated on $K'$ in 
the norm of ${\cal H}_{p,\psi}(K')$ by holomorphic functions from ${\cal H}_{p,\psi}(M')$.
\end{Th}

In the case of coverings of Stein manifolds the results similar to Theorems \ref{te2} and \ref{te3} are proved in [Br4, Theorems 1.8, 1.10].
\sect{\hspace*{-1em}. Proof of Theorem \ref{te1}.}
{\bf 2.1.} We begin the proof with the following auxiliary result.
\begin{Proposition}\label{pr2.1}
For every $z\in M\setminus C_{M}$ and $p\in [1,\infty]$ there is a linear operator
$T_{\psi,z}\in {\cal B}(l_{p,\psi,z}(M'),{\cal H}_{p,\psi}(M'))$ such that
$$
(T_{\psi,z}h)(x)=h(x)\ \ \ {\rm for\ any}\ \ \ h\in l_{p,\psi,z}(M')\ \ \ {\rm and}\ \ \ x\in r^{-1}(z).
$$
\end{Proposition}
(In the sequel we call such $T_{\psi,z}$ a linear interpolation operator.)\\
{\bf Proof.} Let $\widetilde M\subset\subset N$ be a strongly pseudoconvex manifold containing $\overline M$ and $\widetilde M':=r^{-1}(\widetilde M)\subset N'$ be the corresponding covering of $\widetilde M$. Then Theorem 1.1(a) of [Br3] implies that for every function $f\in l_{2,\psi,z}(M')$ there exists $F\in H^{2}_{\psi}(\widetilde M')$ such that $F|_{r^{-1}(z)}=f$. Thus the operator $R_{z}: H^{2}_{\psi}(\widetilde M')\to l_{2,\psi,z}(M')$, $R_{z}g:=g|_{r^{-1}(z)}$, is a linear continuous surjective map of Hilbert spaces. In particular, there is a linear continuous map $S_{\psi,z}:l_{2,\psi,z}(M')\to H^{2}_{\psi}(\widetilde M')$ such that $R_{z}\circ S_{\psi,z}=id$.
\begin{R}\label{re2.2}
{\rm The facts that $R_{z}$ maps $H^{2}_{\psi}(\widetilde M')$ into $l_{2,\psi,z}(M')$ and is continuous easily follow from the uniform continuity of $\log\psi$ and the mean value property for plurisubharmonic functions. Similarly one obtains that the restriction operator $R_{M'}:H^{2}_{\psi}(\widetilde M')\to {\cal H}_{2,\psi}(M')$, $g\mapsto g|_{M'}$, is continuous.}
\end{R}

We set
$$
T_{\psi,z}:=R_{M'}\circ S_{\psi,z}.
$$
Then $T_{\psi,z}$ is the required interpolation operator for $p=2$. Let us prove the result for $p\neq 2$. 

We will naturally identify $r^{-1}(z)$ with $\{z\}\times S$ where $S$ is the fibre of $r$. Let $\{e_{s}\}_{s\in S}$, $e_{s}(z,t)=0$ for $t\neq s$ and $e_{s}(z,s)=(\psi(z,s))^{-1/2}$, be the orthonormal basis of $l_{2,\psi,z}(M')$. We set
$$
h_{s,z}:=T_{\psi,z}(e_{s})\in {\cal H}_{2,\psi}(M').
$$
Then for a sequence $a=\{a_{s}\}_{s\in S}\in l^{2}(S)$ we have
\begin{equation}\label{e2.1}
h_{a}:=\sum_{s\in S}a_{s}h_{s,z}\in {\cal H}_{2,\psi}(M')
\ \ \ {\rm and}\ \ \ |h_{a}|_{2,\psi}\leq c||a||_{l^{2}(S)}.
\end{equation}
We define $F_{s,z}\in {\cal H}_{1,\psi}(M')$ by the formula
\begin{equation}\label{e2.2}
F_{s,z}(w):=\psi(z,s)\ \!h_{s,z}^{2}(w),\ \ \ w\in M'.
\end{equation}
Then (\ref{e2.1}) yields
\begin{equation}\label{e2.3}
\sum_{s\in S}\frac{|F_{s,z}(w)|}{\psi(z,s)}\ \!\psi(w)\leq c^{2},\ \ \ w\in M'.
\end{equation}

Next, for $a=\{a_{s}\}_{s\in S}\in l_{p,\psi,z}(M')$ (i.e.,
$\sum_{s\in S}|a_{s}|^{p}\ \!\psi(z,s):=|a|_{p,\psi,z}^{p}<\infty$) we define
\begin{equation}\label{e2.4}
\widetilde T_{\psi,z}(a):=\sum_{s\in S}a_{s}F_{s,z}.
\end{equation}
For $p=1,\infty$ we set $T_{\psi,z}:=\widetilde T_{\psi,z}$ and show that $T_{\psi,z}$ is the required interpolation operator. In fact, for $p=\infty$ using (\ref{e2.3}) we obtain
$$
\begin{array}{c}
\displaystyle
\sup_{w\in M'}\{|(T_{\psi,z}(a))(w)|\psi(w)\}:=
\sup_{w\in M'}\left\{\left|\sum_{s\in S}a_{s}F_{s,z}(w)\right|\psi(w)\right\}\leq\\
\\
\displaystyle
\left(\sup_{s\in S}\{|a_{s}|\psi(z,s)\}\right)\cdot\left(\sup_{w\in M'}\left\{\sum_{s\in S}\frac{|F_{s,z}(w)|}{\psi(z,s)}\ \!\psi(w)\right\}\right)\leq c^{2}|a|_{\infty,\psi,z}.
\end{array}
$$
Also,
$$
(T_{\psi,z}a)(z,t):=\sum_{s\in S}a_{s}F_{s,z}(z,t):=
a_{t}\psi(z,t)\ \!e_{t}^{2}(z,t)=a_{t}:=a(z,t).
$$
Thus $T_{\psi,z}$ is the required interpolation operator for $p=\infty$.

Similarly for $p=1$ we have from (\ref{e2.1}) with $h_{a}:=h_{s,z}$
$$
\begin{array}{c}
\displaystyle
\sup_{w\in M}\left\{\sum_{y\in r^{-1}(w)}|(T_{\psi,z}(a))(y)|\psi(y)\right\}:=
\sup_{w\in M}\left\{\sum_{y\in r^{-1}(w)}\left|\sum_{s\in S}a_{s}F_{s,z}(y)\right|\psi(y)\right\}\leq\\
\\
\displaystyle
\left(\sum_{s\in S}|a_{s}|\psi(z,s)\right)\cdot\left(\sup_{w\in M}\left\{\sup_{s\in S}\left\{\sum_{y\in r^{-1}(w)}\frac{|F_{s,z}(y)|}{\psi(z,s)}\ \!\psi(y)\right\}\right\}\right)\leq c^{2}|a|_{1,\psi,z}.
\end{array}
$$
Thus $T_{\psi,z}$ is the required interpolation operator for $p=1$. 

Now, using the M. Riesz interpolation theorem (see, e.g., [Ru]) and arguing as in the proof of Lemma 3.2(a) of [Br4] from the cases considered above we obtain that the operator $\widetilde T_{\psi,z}$ maps $l_{p,\psi^{p},z}(M')$ continuously into ${\cal H}_{p,\psi^{p}}(M')$, $1<p<\infty$, and its norm is bounded by $c^{2}$. For such $p$ we define
$$
T_{\psi,z}:=\widetilde T_{\psi^{1/p},z}.
$$
(Note that $\psi^{1/p}$ satisfies the conditions of Theorem \ref{te1}, see section 1.3.)
Then $T_{\psi,z}\in {\cal B}(l_{p,\psi,z}(M'),{\cal H}_{p,\psi}(M'))$ is the required interpolation operator.\ \ \ \ \ $\Box$\\
{\bf 2.2.} We proceed with the proof of the theorem.

Let $E_{0}(M):=M\times {\cal H}_{p,\psi}(M')$ be the trivial holomorphic Banach vector bundle on $M$ with fibre ${\cal H}_{p,\psi}(M')$. Also, by $E_{p,\phi}(M)$ we denote the Banach vector bundle associated with the natural action of $\pi_{1}(M)$ on the fibre $S$ of $r:M'\to M$. For the construction of such a bundle see [Br4, Example 2.2(b)]. Let us recall that the fibre of $E_{p,\phi}(M)$ is the Banach space
$l_{p,\phi}(S)$ of complex functions $f$ on $S$ with norm
$$
||f||_{p,\phi}:=\left(\sum_{s\in S}|f(s)|^{p}\phi(s)\right)^{1/p}
$$
for $p\in [1,\infty)$ and
$$
||f||_{\infty,\phi}:=\sup_{s\in S}\{|f(s)|\phi(s)\}.
$$
Also, $\phi$ is defined by the formula
$$
\phi:=\psi|_{r^{-1}(z_{0})}
$$
for a fixed point $z_{0}\in M$. (Here $\log\psi$ is uniformly continuous on $M'$ with respect to the path metric induced by a Riemannian metric pulled back from $N$.)

It was proved in [Br4, Proposition 2.4] that the Banach space
$B_{p,\phi}(M)$ of bounded holomorphic sections of $E_{p,\phi}(M)$ is naturally isomorphic to ${\cal H}_{p,\psi}(M')$.

Further, for $z\in M$ let $R_{z}$ be the restriction map of functions from ${\cal H}_{p,\psi}(M')$ to the fibre $r^{-1}(z)$.
If we identify ${\cal H}_{p,\psi}(M')$ with the Banach space $B_{p,\phi}(M)$ then $R_{z}$ will be the evaluation map of sections from $B_{p,\phi}(M)$ at $z\in M$. In particular one can define a (holomorphic) homomorphism of bundles $R: E_{0}(M)\to E_{p,\phi}(M)$ which maps $z\times v\in E_{0}(M)$ to the vector $R_{z}(v)$ in the fibre over $z$ of the bundle
$E_{p,\phi}(M)$. Since for every $z\in M\setminus C_{M}$ the space $l_{p,\psi,z}(M')$ is isomorphic to $l_{p,\phi}(S)$, see inequality (2.5) of [Br4], Proposition \ref{pr2.1} implies that
every $R_{z}$ is surjective and, moreover, there exists a linear
continuous map $C_{z}$ of the fibre $E_{p,\phi,z}$ of $E_{p,\phi}(M)$ over $z$ to the fibre $E_{0,z}(M)$ of $E_{0}(M)$ over $z$ such that $R_{z}\circ C_{z}=id$. Repeating literally the arguments of [Br4, section 3.2] we obtain from here that
{\em for every $z\in M\setminus C_{M}$ there is a neighbourhood
$U_{z}\subset\subset M\setminus C_{M}$ of $z$ such that $Ker R|_{U_{z}}$ is biholomorphic to $U_{z}\times Ker R_{z}$ and this biholomorphism is linear on every $Ker R_{x}$ and maps this space onto $x\times Ker R_{z}$, $x\in V_{z}$.} 

The latter shows that the bundle $E_{p,\phi}(M)$ is locally complemented in $E_{0}(M)$ over $M\setminus C_{M}$. Now, if $\Omega\subset M\setminus C_{M}$ is an open Stein manifold, then by the Bungart theorem [B] based on the previous statement one obtains that $E_{p,\phi}(M)$ is complemented in $E_{0}(M)$ over $\Omega$, see [Br4, section 3.2] for details. This means that there is a (holomorphic) homomorphism of bundles
$F:E_{p,\phi}(M)|_{\Omega}\to E_{0}(M)|_{\Omega}$ such that $R\circ F=id$. Moreover, by the definition $F|_{K}$ is bounded over $K\subset\subset\Omega$.

Finally, we set 
$$
L_{z}:=F(z),\ \ \ z\in M.
$$
Then by definition, every $L_{z}$ is a linear continuous map of
$l_{p,\psi,z}(M')$ into ${\cal H}_{p,\psi}(M')$,
the family $\{L_{z}\}$ is holomorphic in $z\in M$, 
$R_{z}\circ L_{z}=id$, and $\sup_{z\in K}||L_{z}||<\infty$.

This completes the proof of the theorem.\ \ \ \ \ $\Box$
\sect{\hspace*{-1em}. Proofs.}
{\bf 3.1. Proof of Theorem \ref{te2}.}
Let $f\in {\cal H}_{p,\psi}(Y')$ be a function where $Y'$ satisfies assumptions of the theorem. These assumptions imply that there is a strongly pseudoconvex manifold $\widetilde M\subset\subset N$ such that $V\subset\subset\widetilde M$. We apply Theorem \ref{te1} with $M'$ substituted for $\widetilde M'$. Then we consider the function
$$
h(z):=L_{z}(f|_{r^{-1}(z)})
,\ \ \ z\in Y.
$$
By [Br4, Proposition 2.4] and by the properties of 
$\{L_{z}\}$ we obtain
that $h$ is a ${\cal H}_{p,\psi}(\widetilde M')$-valued holomorphic function on $Y$. (It can be written as the scalar function of the variables $(z,w)\in Y\times \widetilde M'$.) 
Thus it suffices to prove the extension theorem for the Banach-valued
holomorphic function $h$ on $Y$ extending it to $V$. Evaluating
the extended Banach-valued function at the points $(r(y),y)$, $y\in U'$, we get the required function $F$ (cf. arguments in [Br4, section 4]). Now, the above Banach-valued extension theorem follows directly from the Banach-valued version of the classical Cartan B theorem for Stein manifolds due to Bungart [B].\ \ \ \ \ $\Box$\\
{\bf 3.2. Proof of Theorem \ref{te3}.} We retain the notation of Remark \ref{re3}. By the conditions of the theorem we obtain that 
$X_{K}:=p(K)$ is a holomorphically convex compact subset of
$X_{M}$ that does not contain points $x_{i}$, $1\leq i\leq l$. Then there is a non-degenerate analytic polyhedron $P\subset\subset X_{O}$ containing $K$ and formed by holomorphic functions on $X_{M}$. Now, for $f\in {\cal H}_{p,\psi}(O')$ we consider the function
$$
h(z):=L_{z}(f|_{r^{-1}(z)}),\ \ \ z\in O,
$$
with $\{L_{z}\}$ as in Theorem \ref{te1}. Then
$h$ is a ${\cal H}_{p,\psi}(M')$-valued holomorphic function on $O$. Next, we apply to $h|_{p^{-1}(P)}$ the Weil integral formula (also valid for Banach-valued holomorphic functions). Expanding the kernel in this formula in an analytic series, as in the classical case we obtain that $h$ can be approximated uniformly on $K$ by ${\cal H}_{p,\psi}(M')$-valued holomorphic functions on $M$. Taking restrictions of these functions to the set $\{(r(y),y)\ :\ y\in M'\}$, we obtain the required approximation theorem.\ \ \  \ \ $\Box$

\end{document}